\documentclass[10pt]{article}
\usepackage{latexsym}\usepackage{amsbsy}\usepackage{amssymb}
\usepackage[latin1]{inputenc} 

\newcommand{\ep}{\hspace*{\fill}$\Box$}
\newcommand{\eps}{\varepsilon}
\newcommand{\pr}{\noindent{\bf Proof. }}

\newcommand{\R}{\mathbb R}
\newcommand{\N}{\mathbb N}
\newcommand{\C}{\mathbb C}

\newcommand{\Hbb}{\mathbb H}

\newcommand{\gs}{\ensuremath{{\mathcal G}} }

\newcommand{\es}{\ensuremath{{\mathcal E}} }
\newcommand{\esm}{\ensuremath{{\mathcal E}_M} }

\newcommand{\ns}{\ensuremath{{\mathcal N}} }

\newcommand{\comp}{\subset\subset}
\newcommand{\cinfty}{{\cal C}^\infty}
\newtheorem{thr}{\hspace*{-3mm} \bf}[section]
\newcommand{\bt}{\begin{thr} {\bf Theorem. }}
\newcommand{\et}{\end{thr}}
\newcommand{\bp}{\begin{thr} {\bf Proposition. }}
\newcommand{\bc}{\begin{thr} {\bf Corollary. }}
\newcommand{\blem}{\begin{thr} {\bf Lemma. }}
\newcommand{\bex}{\begin{thr} {\bf Example. }\rm} 
\newcommand{\bexs}{\begin{thr} {\bf Examples. }\rm}
\newcommand{\bd}{\begin{thr} {\bf Definition. }}
\newcommand{\beast}{\begin{eqnarray*}}
\newcommand{\eeast}{\end{eqnarray*}}

\newcommand{\al}{\alpha}

\newcommand{\vphi}{\varphi}

\newcommand{\D}{{\cal D}}

 \newcommand{\G}{{\cal G}}

\newcommand{\beq}{ \begin{equation} }\newcommand{\eeq}{\end{equation} }
\newcommand{\bea}{\begin{eqnarray}}\newcommand{\eea}{\end{eqnarray}}
\newcommand{\beas}{\begin{eqnarray*}}\newcommand{\eeas}{\end{eqnarray*}}
\newcommand{\beqs}{\begin{equation*}}\newcommand{\eeqs}{\end{equation*}}

\newcommand{\bfs}{\boldsymbol}
\newcommand{\ben}{\begin{enumerate}}\newcommand{\een}{\end{enumerate}}
\newcommand{\ba}{\begin{array}}\newcommand{\ea}{\end{array}}
\newcommand{\brem}{\begin{thr} {\bf Remark. }\rm}
\newcommand{\ethi}{\end{thr}}

\newcommand{\supp}{\mathrm{supp}}

\newcommand{\model}{\approx_{\mathcal M}}

\begin{document}

\title{Sheaves of nonlinear generalized functions and manifold-valued distributions}
\author{Michael Kunzinger \footnote{Electronic mail: michael.kunzinger@univie.ac.at}\\ 
        Roland Steinbauer \footnote{Electronic mail: roland.steinbauer@univie.ac.at}\\
         {\small Department of Mathematics, University of Vienna}\\
         {\small Nordbergstr. 15, A-1090 Wien, Austria}\\
        James A.\ Vickers \footnote{Electronic mail: J.A.Vickers@maths.soton.ac.uk}\\
         {\small  University of Southampton, Faculty of Mathematical Studies,}\\ 
         {\small Highfield, Southampton SO17 1BJ, United Kingdom}
       }
\maketitle
                                                                           
\begin{abstract}
This paper is part of an ongoing program to develop a theory of generalized differential 
geometry. We consider the space $\mathcal{G}[X,Y]$ of Colombeau generalized functions 
defined on a manifold $X$ and taking values in a manifold $Y$. 
This space is essential in order to study concepts such as flows of generalized 
vector fields or geodesics of generalized metrics.
We introduce an
embedding of the space of continuous mappings $\mathcal{C}(X,Y)$ into $\mathcal{G}[X,Y]$
and study the sheaf properties of $\mathcal{G}[X,Y]$.
Similar results are obtained for spaces of generalized vector bundle homomorphisms.
Based on these constructions we propose the definition of a
space $\mathcal{D}'[X,Y]$ of distributions on $X$ taking values in $Y$. $\mathcal{D}'[X,Y]$ is
realized as a quotient of a certain subspace of $\mathcal{G}[X,Y]$.

\medskip

\medskip 

\noindent{\footnotesize This work was supported by project P16742 and START-project Y-237 of the 
Austrian Science Fund}

\end{abstract}

\section{Introduction}\label{intro}

Non-linear generalized functions in the sense of J.F.\ Colombeau 
(\cite{c1,c2,Cbull}) are an extension of the (linear) theory of 
distributions providing maximal consistency with respect to classical
analysis in light of L. Schwartz's impossibility result
(\cite{Schw}). While originally used as a tool for studying non-linear partial differential
equations (see \cite{MOBook} for a survey), applications of a more geometric
nature, in particular, in Lie group analysis of differential equations
(e.g.\ \cite{symm, book}) and in general relativity (see \cite{SV} for a recent review)
have led to the development of a {\em geometric} theory of non-linear generalized functions
(see \cite{book} for an overview).

According to general relativity 
the curvature of space-time is given by Einstein's equations which form a non-linear system of second order partial differential equations for the metric. In order to have a well defined space-time one therefore requires sufficient differentiability of the metric for Einstein's equations to make sense. However there are a number of physically important solutions for which the metric does not possess the necessary level of differentiability. Examples of this include space-times with pp-waves and cosmic strings. In order to enlarge the class of space-times one can deal with, so that they include those of physical importance, it is necessary to develop a theory of distributional geometry. However because the curvature is a non-linear function of the derivatives of the components of the metric such a theory has to be produced using the theory of non-linear generalized functions. A key step in this program was the construction of a diffeomorphism invariant scalar theory in the so called `full' setting where one has a canonical embedding of distributions into the algebra \cite{found,vim,Jel2,Jel3}. This built upon the pioneering work of \cite{AB,CM,Jel}. However in order to develop a theory of generalized differential geometry one needs to go beyond this and have a description of generalized tensor fields, a topic of ongoing research.
On the other hand for the so-called `special' setting (which will be the framework
of this article) building on \cite{RD} a theory of generalized sections in vector bundles 
was introduced in \cite{ndg}. It extends the distributional 
theory of De Rham and Marsden (\cite{deR,marsden}) and was used to introduce a generalized 
(pseudo-)Riemannian geometry in \cite{gprg}. Moreover, the need to consider geodesics of a generalized (pseudo-)\-Rie\-man\-nian metric or 
the flow of a generalized vector field made it essential to consider generalized 
functions taking {\em values} in a manifold (a concept not available in classical distribution theory). 
To this end, in \cite{gfvm} the space $\G[X,Y]$ 
of generalized functions defined on a manifold $X$ and taking values in a manifold $Y$ 
was introduced as well as the space $\mathrm{Hom}_{\G}(E,F)$ of generalized vector bundle homomorphisms.
The notions of \cite{gfvm} were used in the description of geodesics in impulsive gravitational pp-waves
of general relativity (see \cite{book}, Chapter 5 for an overview) as well as in the study of flows of
generalized vector fields (\cite{flows}). On the other hand in \cite {gfvm2} the construction of \cite{gfvm}
was completed and turned into a functorial theory. In particular, several global characterizations of the
spaces $\G[X,Y]$ and $\mathrm{Hom}_{\G}(E,F)$ were given. 
In this work we develop this line
of research further and consider sheaf properties of these spaces as well as the question of the embedding 
of `regular` (e.g. continuous resp.\ locally bounded) functions into the respective spaces.
Based on this we propose the definition of a space of manifold-valued distributions.
Our approach extends the sheaf-theoretic study of real-valued
Colombeau generalized functions given in \cite{RD, DP, DP2, marti}, cf.\ also
\cite{DHPV} for an alternative setting.

In some more detail, in section \ref{sheaves} we prove that  $\G[\,\_,Y]$ and $\mathrm{Hom}_{\G}(\,\_,F)$ 
are sheaves of sets. Section \ref{embeddings} is devoted to proving that the space ${\mathcal C}(X,Y)$ of 
continuous functions from $X$ to $Y$ is embedded into $\G[X,Y]$ as well as the analogous 
statement for vector bundle homomorphisms. Furthermore we discuss the (im)possibility of embedding
locally bounded functions into $\G[X,Y]$ and provide an embedding for vector bundle homomorphisms 
which display a `more singular' behavior in their fiber component. Finally, in section \ref{D'(X,Y)} 
we propose the introduction of a space $\D'[X,Y]$ of `distributions' taking values in $Y$---a
notion which does not exist classically. Inspired by the sequential approach to distribution theory on Euclidean space (see \cite{Mik,AMS,Temp}) we utilize our constructions to define a quotient of a certain
subspace of $\G[X,Y]$ which serves as an analog to these distributional spaces. We conclude this work
by investigating the basic properties of $\D'[X,Y]$.

In the remainder of this introduction we recall some notions from \cite{gfvm,gfvm2} which will be needed
in the sequel and fix some notation. Our main reference on non-linear generalized functions is \cite{book}.

Throughout this paper $X$ and $Y$ denote smooth Hausdorff manifolds with countable basis of 
dimension $m$ and $n$, respectively. We set $I=(0,1]$ and 
${\mathcal E}(X):=\cinfty(I\times X,\C)$. Elements of ${\mathcal E}(X)$ will as usual
be denoted in the form $(u_\eps)_\eps$ to emphasize the role of $\eps$ as a
(regularization) parameter. Note, however, that here and in what follows we 
suppose smooth dependence on $\eps$. This additional assumption leaves unchanged all
properties of the spaces of generalized functions as described in \cite{book, gfvm, gfvm2},
yet it will be crucial for the constructions to follow. Independent reasons for
this choice of basic space are certain algebraic simplifications: for
polynomials with generalized coefficients to only have classical solutions one needs
at least continuous dependence on $\eps$, cf. \cite{MOBook}, Prop. 12.2.\ or
\cite{Biag}, Prop.\ 1.10.5. A similar statement holds for solutions of
polynomial ordinary differential equations, see \cite{Biag}, Cor.\ 1.10.9.

The special Colombeau algebra of generalized functions $\gs(X)$ on $X$ is defined 
as the quotient $\esm(X)/\ns(X)$ of moderate modulo negligible nets in ${\mathcal E}(X)$. 
Here the latter notions are defined by (denoting by ${\cal P}(X)$ the space of linear 
differential operators on $X$)
\beas
        \esm(X)&:=&\{ (u_\eps)_\eps\in{\mathcal E}(X):\ 
        \forall K\subset\subset X,\ \forall P\in{\cal P}(X)\ \exists N\in\N:
  \\&&\hphantom{(u_\eps)_\eps\in{\mathcal E}(X):\ \forall K\subset\subset X,\ \forall 
     P\in{\cal P}(\}} 
        \sup_{p\in K}|Pu_\eps(p)|=O(\eps^{-N})\}
   \\
        \ns(X)&:=& \{ (u_\eps)_\eps\in\esm(X):\ 
        \forall K\subset\subset X,\ \forall q \in\N_0:\
        \sup_{p\in K}|u_\eps(p)|=O(\eps^{q}))\}\,.
\eeas
Since we are going to work entirely in the `special' setting of 
Colombeau's construction we omit this term henceforth.
$\gs(\_)$ is a fine sheaf of differential algebras
with respect to the Lie derivative along smooth vector fields defined by 
$L_\xi u=[(L_\xi u_\eps)_\eps]$ (\cite{RD, ndg}). The construction can be 
appropriately localized, that is $u$ is in $\gs(X)$ if and only 
if $u\circ\psi_\al\in\gs(\psi_\al(V_\al))$ (the local Colombeau algebra on 
$\psi_\al(V_\al)$) for all charts $(V_\al,\psi_\al)$. 
$\cinfty(X)$ is a sub-algebra of $\gs(X)$ and there exist injective sheaf morphisms 
embedding $\D'(X)$, the space of Schwartz distributions on $X$, into $\gs(X)$.

The space $\gs[X,Y]$ of compactly bounded (c-bounded) generalized Co\-lom\-beau functions on $X$ 
taking values in $Y$ is defined similarly as a quotient of the set $\es_M[X,Y]$ 
of moderate, c-bounded maps from $X$ to $Y$ by a certain equivalence 
relation. However, in the absence of a linear structure in the target space the definition of
the respective asymptotics becomes more involved.
\bd 
\begin{itemize}
\item [(a)]
 $\esm[X,Y]$ is defined as the set of all $(u_\eps)_\eps \in \cinfty(I\times X,Y)$ 
 satisfying
 \begin{itemize}
  \item[(i)] $\forall K\comp X\ \exists \eps_0>0\  \exists K'\comp Y \  \forall
             \eps<\eps_0:\ u_\eps(K) \subseteq K'$ (c-boundedness).
  \item[(ii)]
   $\forall k\in\N$, for each chart $(V,\vphi)$ in $X$, each 
   chart $(W,\psi)$ in $Y$, each $L\comp V$ and each $L'\comp W$
   there exists $N\in \N$  with
   $$\sup\limits_{x\in L\cap u_\eps^{-1}(L')} \|D^{(k)}
   (\psi\circ u_\eps \circ \vphi^{-1})(\vphi(x))\| =O(\eps^{-N}).
   $$
 \end{itemize}
 \item[(b)] $(u_\eps)_\eps$ and $(v_\eps)_\eps \in \esm[X,Y]$ are
  called equivalent, $(u_\eps)_\eps \sim (v_\eps)_\eps$,
  if the following conditions are satisfied:
  \begin{itemize}
    \item[(i)]  For all $K\comp X$, $\sup_{p\in K}d_h(u_\eps(p),v_\eps(p)) \to 0$
    ($\eps\to 0$)
    for some (hence every) Riemannian metric $h$ on $Y$.
    \item[(ii)] $\forall k\in \N_0\ \forall m\in \N$,
    for each chart
    $(V,\vphi)$
    in $X$, each chart $(W,\psi)$ in $Y$, each $L\comp V$
    and each $L'\comp W$:
    $$
     \sup\limits_{x\in L\cap u_\eps^{-1}(L')\cap v_\eps^{-1}(L')}\!\!\!\!\!\!\!\!\!\!\!\!\!\!\!\!\!\!\!
     \|D^{(k)}(\psi\circ u_\eps\circ \vphi^{-1}
     - \psi\circ v_\eps\circ \vphi^{-1})(\vphi(x))\|
     =O(\eps^m).
    $$
   \end{itemize}
  \item[(c)]  The space of generalized functions from $X$ to $Y$ is defined as
   \[\gs[X,Y]\,:=\,\es_M[X,Y]/\sim\ .\]
 \end{itemize}
\et

The following characterization result has been established in \cite{gfvm2}, Prop.\ 3.2 
and Th.\ 3.3 and will be repeatedly used throughout this work.
\bt \label{mainchar}
\begin{itemize}
\item[(i)] Let $(u_\eps)_\eps \in \cinfty(I\times X, Y)$. Then $(u_\eps)_\eps \in \esm[X,Y]$
if and only if $(f\circ u_\eps)_\eps \in \esm(X)$ for all $f\in \cinfty(Y)$.
\item[(ii)] If $(u_\eps)_\eps$, $(v_\eps)_\eps \in \esm[X,Y]$, then 
$(u_\eps)_\eps \sim (v_\eps)_\eps$ if and only if $(f\circ u_\eps - f\circ v_\eps)_\eps \in 
\ns(X)$ for all $f\in \cinfty(Y)$.
\end{itemize}
\et
 
Finally, we turn to the definition of generalized vector bundle homomorphisms
(e.g., tangent maps of manifold valued generalized functions).
Vector bundles with base space $X$ will be denoted 
$(E,X,\pi_X)$. A vector bundle chart $(V,\Phi)$ over a chart  $(V,\varphi)$ of $X$ will be written in the form 
$\Phi(e) = (\vphi(p),\bfs{\vphi}(e)) \in \vphi(V)\times \R^{n'}$
where $p=\pi_X(e)$. The space of smooth vector bundle homomorphisms 
from $E$ to $(F,Y,\pi_Y)$ will be called  $\mathrm{Hom}(E,F)$.
If $f\in \mathrm{Hom}(E$, $F)$ we write $\underline{f}: X\to Y$ for the smooth 
map induced on bases, i.e., $\pi_Y\circ f = \underline{f}\circ \pi_X$.
Local vector bundle homomorphisms with respect to 
vector bundle charts $(V,\Phi)$ of $E$ and $(W,\Psi)$ of $F$, i.e.,
$f_{\mathrm{\Psi}\mathrm{\Phi}}:=
\mathrm{\Psi}\circ f \circ \mathrm{\Phi}^{-1}: \vphi(V\cap \underline{f}^{-1}(W))
\times \R^{m'} \to \psi(W) \times \R^{n'}$
will be written in the form
\begin{equation}\label{vbhomloc}
f_{\mathrm{\Psi}\mathrm{\Phi}}(x,\xi) = 
(f_{\mathrm{\Psi}\mathrm{\Phi}}^{(1)}(x),f_{\mathrm{\Psi}\mathrm{\Phi}}^{(2)}(x)
\cdot\xi)\,.
\end{equation}
Here, $f_{\mathrm{\Psi}\mathrm{\Phi}}^{(1)} = \underline{f}_{\psi\vphi} := \psi\circ \underline f \circ\vphi^{-1}$. Finally we are ready to define the space of generalized vector bundle homomorphisms.

\bd\label{longdef}
\begin{itemize}
\item[(a)]
${\esm}^{\mathrm{VB}}[E,F]$ is the set of all $(u_\eps)_\eps$ 
$\in$ $\mathrm{Hom}(E,F)^I$ depending smoothly on $\eps$ and satisfying
\begin{itemize}
\item[(i)] $(\underline{u_\eps})_\eps \in \esm[X,Y]$.
\item[(ii)] $\forall k\in \N_0\
\forall (V,\Phi)$
vector bundle chart in $E$,
$\forall (W,\Psi)$ vector bundle chart in $F$,
$\forall L\comp V\
\forall L'\comp W\ \exists N\in \N\ \exists \eps_1>0\
\exists C>0$ with
$$
\|D^{(k)}
(u_{\eps \mathrm{\Psi}\mathrm{\Phi}}^{(2)}(\vphi(p)))\|
\le C\eps^{-N}
$$
for all $\eps<\eps_1$ and all $p\in L\cap\underline{u_\eps}^{-1}(L')$, where $\|\,.\,\|$ 
denotes any matrix norm.
\end{itemize}
\item[(b)]\label{homgequ}
$(u_\eps)_\eps$, $(v_\eps)_\eps \in {\esm}^{\mathrm{VB}}[E,F]$
are called $vb$-equivalent, $((u_\eps)_\eps \sim_{vb} (v_\eps)_\eps)$
if
\begin{itemize}
\item[(i)] $(\underline{u_\eps})_\eps \sim (\underline{v_\eps})_\eps$ in
$\esm[X,Y]$.
\item[(ii)] $\forall k\in \N_0\ \forall m\in \N\ \forall (V,\Phi)$
vector bundle chart in
$E$, $\forall (W,\Psi)$ vector bundle chart in $F$,
$\forall L\comp V\ \forall L'\comp W
\ \exists \eps_1>0\ \exists C>0$ such that:
$$
\|D^{(k)}(u_{\eps \mathrm{\Psi}\mathrm{\Phi}}^{(2)}
-v_{\eps \mathrm{\Psi}\mathrm{\Phi}}^{(2)})(\vphi(p))\|
\le C\eps^{m}
$$
for all $\eps<\eps_1$ and all $p\in L\cap\underline{u_\eps}^{-1}(L')
\cap\underline{v_\eps}^{-1}(L')$.
\end{itemize}
\item[(c)] The space of generalized vector bundle homomorphisms is defined by
 $$\mathrm{Hom}_{\gs}[E,F] := {\esm}^{\mathrm{VB}}[E,F]\big/\sim_{vb}.$$
\end{itemize}
\et
 
For $u\in \mathrm{Hom}_{\gs}[E,F]$, $\underline{u} :=[(\underline{u}_\eps)_\eps]$
is a well-defined element of $\gs[X,Y]$ uniquely characterized by $\underline{u}
\circ\pi_X = \pi_Y\circ u$. The tangent map $Tu:=[(Tu_\eps)_\eps]$ of any 
$u\in\gs[X,Y]$ is a well-defined element of $\mathrm{Hom}_{\gs}[TX,TY]$.

\section{The sheaf property of $\gs[X,Y]$ and $\mathrm{Hom}_\gs(E,F)$}\label{sheaves}
Our aim in this section is to establish that $\gs[\,\_\,,Y]$ is a sheaf of sets. 
Compared to the case of the Colombeau algebra $\gs(X)$ (where the sheaf property
can basically be derived by lifting the local result (\cite{c2}, \S 1.3) to the manifold
(\cite{RD, ndg})) the main obstacle in the present setting is the lack of algebraic
structure on the target space $Y$. Given an open cover $\{U_\alpha \mid \alpha \in A\}$
of $X$ it follows directly from the definition (or also from Th.\ \ref{mainchar} (ii)) that
if $u$, $v \in \gs[X,Y]$ and $u|_{U_\alpha} = v|_{U_\alpha}$ for all $\alpha\in A$, then 
$u=v$. The second defining property: 
$$
\left.
\mbox{
\begin{minipage}{100mm}
Given a family $\{u_\alpha\in \gs[U_\alpha,Y] \mid \alpha \in A\}$
such that $u_\alpha|_{U_\alpha\cap U_\beta} = u_\beta|_{U_\alpha\cap U_\beta}$ for all 
$\alpha$, $\beta$ with $U_\alpha\cap U_\beta \not=\emptyset$ then there exists some
$u\in \gs[X,Y]$ with $u|_{U_\alpha} = u_\alpha$ for all $\alpha\in A$.
\end{minipage}
}
\right \}(*)
$$
however, cannot be established similar to the algebra-setting: the tools for
glu\-ing together locally defined elements of our quotient spaces (e.g., partitions of unity) are absent in
the manifold-valued case. Our strategy therefore will be to first embed the target
manifold in some Euclidean space by a Whitney embedding, do the gluing via partitions
of unity in the surrounding vector space and then project back onto the target manifold $Y$
using the retraction map of a tubular neighborhood of $Y$. To this end we first provide an
alternative characterization of $\gs[X,Y]$ in case $Y$ is a submanifold of some $\R^s$.

\bd \label{gtildedef} Let $Y$ be a submanifold of $\R^s$. We define
$\tilde \gs[X,Y]$ to be the subset of $\gs(X)^s$ consisting of those 
$u\in \gs(X)^s$ which possess a representative  $(u_\eps)_\eps$ satisfying:
\begin{itemize}
\item[(i)] $u_\eps(X) \subseteq Y$ for all $\eps\in I$.
\item[(ii)] For each $K\comp X$ there exist $\eps_0 > 0$ and $K'\comp Y$ such that
$u_\eps(K)\subseteq K'$ for all $\eps<\eps_0$. 
\end{itemize}
\et
\bp \label{gtildeprop}
Let $X$, $Y$ be smooth manifolds and let $i:Y\hookrightarrow \R^s$ be an 
embedding of $Y$. Then the push-forward $i_*: \gs[X,Y] \to \tilde \gs[X,i(Y)]$, 
$i_*(u) = i\circ u$ is a bijection which commutes with restrictions to open sets.  
In particular, if $Y$ is a submanifold of some
$\R^s$ we may identify $\gs[X,Y]$ with $\tilde \gs[X,Y]$.
\et
\pr We first note that $i_*$ is well-defined: for $(u_\eps)_\eps\in \esm[X,Y]$ 
Th.\ \ref{mainchar} (i) implies that $(i\circ u_\eps)_\eps \in  \esm(X)^s$. Also, properties (i) and
(ii) of Def.\ \ref{gtildedef} follow immediately. Suppose now that $(u_\eps)_\eps \sim (v_\eps)_\eps$
for $(u_\eps)_\eps, \, (v_\eps)_\eps \in \esm[X,Y]$. Then by Th.\ \ref{mainchar} (ii) we have 
that $(i_j\circ u_\eps - i_j\circ v_\eps)_\eps \in \ns(X)$ for $1\le j\le s$, so 
$[(i\circ u_\eps)_\eps] = [(i\circ v_\eps)_\eps]$ in $\gs(X)^s$. Moreover, it follows
directly from the definitions that $i_*$ commutes with restrictions to open sets.

$i_*$ is injective: suppose that $i_*([(u_\eps)_\eps]) = i_*([(v_\eps)_\eps])$, i.e., 
$(i\circ u_\eps)_\eps - (i\circ v_\eps)_\eps \in \ns(X)^s$. An application of the mean 
value theorem shows that this entails  $(g\circ i\circ u_\eps)_\eps - (g\circ i\circ v_\eps)_\eps \in \ns(X)$
for all $g\in \cinfty(i(Y))$. Since $i: Y \to i(Y)$ is
a diffeomorphism, it follows that any $g\in \cinfty(i(Y))$ is of the form $f\circ i^{-1}$
for some $f\in \cinfty(Y)$. Hence Th.\ \ref{mainchar} (ii) gives $(u_\eps)_\eps \sim (v_\eps)_\eps$.

$i_*$ is surjective: Let $\tilde u \in \tilde \gs[X,i(Y)]$ with representative $(\tilde u_\eps)_\eps$
satisfying (i) and (ii) of Def.\ \ref{gtildedef}. Then $u_\eps:= i^{-1}\circ \tilde u_\eps$ defines
(by Th. \ref{mainchar} (i)) an element of $\esm[X,Y]$ whose image under $i_*$ is $\tilde u$. \ep

After these preparations we can now prove:

\bt \label{sheafth} $\gs[\,\_\,,Y]$ is a sheaf of sets.
\et 
\pr By Whitney's embedding theorem (cf.\ \cite{hirsch} or \cite{lee}) there exists an
embedding of $Y$ into some $\R^s$.
Due to Prop.\ \ref{gtildeprop} and our preparatory statements at the beginning of this section
it therefore
suffices to suppose that $Y$ is in fact a submanifold of some $\R^s$ and to establish property
$(*)$ for $\tilde \gs[X,Y]$. Thus let $\{U_\alpha \mid \alpha\in A\}$ be an open cover of $X$
and let $u^\alpha \in \tilde \gs[U_\alpha,Y]$ such that $u^\alpha|_{U_\alpha\cap U_\beta} = 
u^\beta|_{U_\alpha\cap U_\beta}$ for all $\alpha$, $\beta$ with $U_\alpha\cap U_\beta \not=\emptyset$.
Since $X$ is Hausdorff and second countable, it is $\sigma$-compact, and, in particular, paracompact
and Lindel\"of. Without loss of generality we may therefore suppose that $A$ is countable and
that $\{U_\alpha \mid \alpha \in A\}$ is locally finite. Let $\{K_l\mid l\in \N\}$ be an exhaustive
sequence of compact sets in $X$ with $K_l \subseteq K_{l+1}^\circ$ for all $l$. Again without loss
of generality we may suppose that $\{U_\alpha \mid \alpha \in A\}$ is a refinement of 
$\{K_l^\circ\mid l\in \N\}$.

Let $T$ be an open tubular neighborhood of $Y$ in $\R^s$ (again see \cite{hirsch} or \cite{lee}) 
and denote by $r: T\to Y$, $r|_Y = \mathop{id}_Y$ the corresponding retraction. Choose a closed
tubular neighborhood $T'\subseteq T$ and a smooth map $\tilde r: \R^s\to \R^s$ such that 
$\tilde r|_{T'} = r$. Let $\{\chi_\alpha \mid \alpha\in  A\}$ be a partition of unity with 
$\mathrm{supp} \chi_\alpha \comp U_\alpha$ for each $\alpha\in A$. For $\eps\in I$ we set
$$
w_\eps := \tilde r\circ \left( \sum_{\alpha\in A} \chi_\alpha u_\eps^\alpha\right)\,.
$$
Then $(w_\eps)_\eps \in \esm[X,\R^s]$. Fix $l\in \N$ and let $\alpha_1,\dots,\alpha_k$
be the finitely many indices with $\mathrm{supp}\chi_{\alpha_i} \cap K_l \not=\emptyset$
($1\le i\le k$). Choose $\eps_l'>0$ and $K_l'\comp Y$ such that $u_\eps^{\alpha_i}(\supp \chi_{\alpha_i}
\cap K_l) \subseteq K_l'$ for all $1\le i\le k$ and all $\eps<\eps_l'$. For each $y\in Y$ choose 
some $R_y>0$ such that the ball $B_{R_y}(y)$ of radius $R_y$ around $y$ in $\R^s$ is contained in
$T'$. Since $K_l'$ is compact there exists some $\delta>0$ (the Lebesgue number of the covering 
$\{B_{R_y}(y)\cap Y \mid y\in K_l'\}$ of $K_l'$) such that any subset of $K_l'$ with diameter 
less than $\delta$ lies entirely within one $B_{R_y}(y)$.

Since $u^{\alpha_i} - u^{\alpha_j}$ is negligible on $U_{\alpha_i}\cap U_{\alpha_j}$ (in case this set
is nonempty) we may choose some $\eps_l < \eps_l'$ such that $|u_\eps^{\alpha_i}(x) - u_\eps^{\alpha_j}(x)| <\delta$
whenever $x\in K_l\cap \supp \chi_{\alpha_i}\cap \supp \chi_{\alpha_j}$ and $\eps< \eps_l$
($1\le i\le k$). Hence (using the convexity of $B_{R_y}(y)$) for each $x\in K_l$ and each 
$\eps<\eps_l$ there exists some $y\in K_l'$ such that 
$$
\sum_{\alpha\in A} \chi_\alpha(x) u_\eps^\alpha(x) \in B_{R_y}(y) \subseteq T'\,.
$$
Therefore,
$$
w_\eps(x) =  r\circ \left( \sum_{\alpha\in A} \chi_\alpha(x) u_\eps^\alpha(x)\right) 
\qquad \forall x\in K_l \ \forall \eps < \eps_l\,.
$$

If $\beta\in A$ is such that $U_\beta \subseteq K_l^\circ$ and $\eps<\eps_l$ then for each $L\comp U_\beta$
and each $x\in L$ we have
\begin{eqnarray*}
|w_\eps(x) - u_\eps^\beta(x)| &=& \left|\, r\circ \left( \sum_{\alpha\in A} \chi_\alpha(x) u_\eps^\alpha(x)\right)
-r(u_\eps^\beta(x))  \right| \le \\
&\le& \|D\tilde r \|_{L^\infty(\mathrm{ch}(K_l'))} \left|  \sum_{\alpha\in A} \chi_\alpha(x) u_\eps^\alpha(x) - u_\eps^\beta(x) \right|\,,
\end{eqnarray*}
where $\mathrm{ch}(K_l')$, the convex hull of $K_l'$ is itself compact. Since on $L$ the last factor in this estimate
vanishes faster than any power of $\eps$ we have $w|_{U_\beta} = u^\beta$
for all $\beta$ with $U_\beta\subseteq K_l^\circ$.

Choose a smooth function $\eta: X \to \R$ such that $0<\eta(x)\le \eps_l$ for all $x\in K_l\setminus K_l^\circ$
($K_0:=\emptyset$) (see, e.g., \cite{book}, Lemma 2.7.3). Moreover, let $\nu: \R^+ \to [0,1]$ be a smooth function
satisfying $\nu(x)\le x$ for all $x$ and
$$
\nu(x) = \left\{
\begin{array}{ll}
x & 0 \le x \le \frac{1}{2} \\
1 & x \ge \frac{3}{2}                     
\end{array}
\right.
$$
For $(\eps,x)\in I\times X$ we set $\mu(\eps,x) := \eta(x) \nu\left(\frac{\eps}{\eta(x)}\right)$. Finally, we set
$u_\eps(x) := w_{\mu(\eps,x)}(x)$ for $(\eps,x)\in I\times X$. Then $(\eps,x) \mapsto u_\eps(x) \in \cinfty(I\times X,Y)$
(it is here that we need smooth dependence of representatives on $\eps$). Furthermore,
$$
u_\eps|_{K_l^\circ} = w_\eps|_{K_l^\circ} \quad \mbox{ for } 
\eps \le \frac{1}{2} \min_{x\in K_l}\eta(x)\,, \ l\in \N\,.
$$ 
Therefore,
$u = [(u_\eps)_\eps] \in \tilde\gs[X,Y]$ and $u|_{U_\beta} = u^\beta$ for all $\beta\in A$. 
\ep

\brem \label{mgrem}
The method of gluing via the function $\mu$ in the above proof also allows us to establish
the equality of the space of generalized functions taking values in an open subset of $\R^n$
as introduced in \cite{AB} (with smooth dependence on $\eps$) with our setting. 
In fact, let $\Omega\subseteq \R^m$ and 
$\Omega' \subseteq \R^n$ be open sets and let $u\in \gs(\Omega)^n$ such that $u$ possesses
a representative $(u_\eps)_\eps$ satisfying: $\forall K\comp \Omega$ $\exists K'\comp \Omega'$
$\exists \eps_0>0$ such that $u_\eps(K) \subseteq K'$ for all $\eps < \eps_0$. Choose an 
exhaustive sequence $\{K_l\mid l\in \N\}$ of $\Omega$ and corresponding $K_l' \comp \Omega'$
and $\eps_l$ as above. Then defining $\mu$ as in the proof of Th.\ \ref{sheafth}, $\tilde u_\eps(x):=
u_{\mu(\eps,x)}(x)$ defines a representative of $u$ such that $u_\eps(\Omega) \subseteq \Omega'$ 
for all $\eps$. This shows that $\gs_*(\Omega,\Omega')$ in the sense of \cite{AB} can be identified
with $\gs[\Omega,\Omega']$.
\et

In what follows, we want to utilize Th.\ \ref{sheafth} to establish the sheaf property of the space
$\mathrm{Hom}_\gs(E,F)$ of generalized vector bundle homomorphisms. To this end we need some preparatory
constructions for smooth vector bundle homomorphisms.

Let $f: X\to Y$ be any smooth map. Then $f$ can naturally be extended to a vector
bundle homomorphism $\bar f \in \mathrm{Hom}(E,F)$ by defining its action on the fibers
of $E$ to be $0$ (i.e., any local representative $\bar f_{\Psi\Phi}$ of $\bar f$ is of the 
form $(x,v) \mapsto (f_{\psi\vphi}(x),0)$). Suppose now that $U$, $V$ are open subsets of $X$ 
with $V \subset \bar V \comp U$ and let $u \in \mathrm{Hom}(E|_U,F)$ be such that $\underline{u}=
f|_U$. Choose a bump function $\chi\in \cinfty(X)$ such that $\supp \chi \comp U$ and $\chi|_V
\equiv 1$. 
Then
$\chi\cdot u:= e \mapsto \chi(\pi_X(e))u(e)$ (fiber-wise product) defines an element of $\mathrm{Hom}(E|_U,F)$.
Moreover, there is a unique element $v=\chi \bullet_f u$ of $\mathrm{Hom}(E,F)$ such that $v|_U = 
\chi\cdot u$ and $v|_{X\setminus U} = \bar f|_{X\setminus U}$. 
Then $v|_V = u|_V$. We will use these notations in the proof of the following result.
\bt \label{homsheaf} $\mathrm{Hom}_\gs(\pi_X^{-1}(\,\_\,),F)$ is a sheaf of sets on $X$.
\et
\pr As in the case of manifold-valued generalized functions, it follows directly from the
definitions that if $\{U_\alpha\mid \alpha \in A\}$ is an open cover and $u,\, v \in \mathrm{Hom}_\gs(E,F)$
are such that  $u|_{U_\alpha} = v|_{U_\alpha}$ for all $\alpha\in A$ then $u=v$ (Here and in what follows
we abbreviate $u|_{\pi_X^{-1}(U_\alpha)}$ by $u|_{U_\alpha}$).

Suppose that $\{u_\alpha \in \mathrm{Hom}(E|_{U_\alpha}) \mid \alpha\in A\}$ forms a coherent family, i.e., 
$u_\alpha|_{U_\alpha\cap U_\beta} = u_\beta|_{U_\alpha\cap U_\beta}$ for all 
$\alpha$, $\beta$ with $U_\alpha\cap U_\beta \not=\emptyset$. Then $\{\underline u_\alpha \mid \alpha \in A\}$
forms a coherent family in $\gs[X,Y]$, so by Th.\ \ref{sheafth} there exists a unique element $w \in \gs[X,Y]$
such that $w|_{U_\alpha} = \underline u_\alpha$ for all $\alpha\in A$. Thus each $u_\alpha$ is an
element of $\mathrm{Hom}_w(E|_{U_\alpha},F)$, the space of generalized vector bundle homomorphisms with
base component $w$ (cf.\ \cite{gfvm2}, Sec.\ 5). Then by \cite{gfvm2}, Prop.\ 5.7, for each $\alpha \in A$
we may choose a representative $(u^\alpha_\eps)_\eps$ of $u_\alpha$ such that $\underline u^\alpha_\eps = 
w_\eps|_{U_\alpha}$ for all $\eps \in I$. (We note that in order to adapt the proof of Prop.\ 5.7 in
\cite{gfvm2} to the present setting of smooth $\eps$-dependence, a `gluing function' $\mu$ as in 
the proof of Th.\ \ref{sheafth} has to be employed).

Choose now a partition of unity $\{\chi_j\mid j\in \N\}$ subordinate to $\{U_\alpha\mid \alpha\in A\}$
such that $\supp \chi_j \comp U_{\alpha_j}$ for each $j$. For each $\eps\in I$ we define $u_\eps$
as the following (locally finite) fiber-wise sum:
$$
u_\eps:= \sum_{j\in \N} \chi_j \bullet_{w_\eps} u^{\alpha_j}_\eps\,.
$$
Then $u=[(u_\eps)_\eps] \in \mathrm{Hom}_\gs(E,F)$ and $\underline{u} = w$. By \cite{gfvm2},
Th.\ 4.2, in order to show that $u|_{U_\alpha} = u_\alpha$ it suffices to establish
Def.\ \ref{longdef} (b) (ii) for $k=0$ (i.e., we do not have to take into account derivatives).
This, however, is immediate from the coherence of the family $\{u_\alpha\mid \alpha\in A\}$
and the fact that $\{\chi_j\mid j\in \N\}$ is a partition of unity. \ep

\section{Embeddings} \label{embedding}\label{embeddings}
Our aim in this section is to construct embeddings of spaces of 
continuous (resp.\ even more singular) mappings into spaces of
manifold-valued generalized functions. The basic idea (similar to
a procedure introduced in \cite{N}, Part A) is to employ a Whitney embedding of
the target space into some $\R^s$ and then use convolution for
smoothing. The retraction map of a tubular neighborhood of $Y$
in $\R^s$ will then be used to project the resulting nets of smooth
functions back to $Y$.

As was already pointed out in \cite{gfvm}, there is a canonical embedding $\sigma$ of $\cinfty(X,Y)$ into
$\gs[X,Y]$, $\sigma: u \to [(u)_\eps]$. The following result extends this embedding to the space
of continuous mappings from $X$ to $Y$.

\bt \label{embth1}
There exists an embedding $\iota: {\mathcal C}(X,Y) \hookrightarrow \gs[X,Y]$ with the
following properties:
\begin{itemize}
\item[(i)] $\iota$ is a sheaf morphism.
\item[(ii)] $\iota|_{\cinfty(X,Y)} = \sigma$.
\item[(iii)] $\iota(u)_\eps$ converges to $u$ uniformly on compact sets.
\end{itemize}
\et
\pr Using a Whitney embedding of $Y$ in some $\R^s$, Prop.\ \ref{gtildeprop} and the above remarks
show that without loss of generality we may suppose that $Y$ is a submanifold of $\R^s$ and that
it suffices to embed $\mathcal{C}(X,Y)$ into $\tilde \gs[X,Y]$. 

By \cite{ndg}, Th.\ 1 (or \cite{book}, Th.\ 3.2.10, see also \cite{RD}, \S 15 for an equivalent construction
based on de Rham regularizations) there exists an injective sheaf morphism $\tilde \iota: \D'(X,\R^s)
\hookrightarrow \gs(M)^s$. This embedding is based on regularization via convolution with a mollifier
in charts of a given atlas, patched together through a partition of unity. In this way, the convergence
properties of the respective regularizations of 
continuous mappings are preserved. 
In particular, $\tilde \iota(u)_\eps$ converges 
uniformly on compact sets 
to $u$ 
for $u \in {\mathcal C}(X,\R^s)$.

In what follows we use the notations of the proof of Th.\ 
\ref{sheafth}. Let $u\in {\mathcal C}(X,Y)$.  For each $l\in \N$ we
choose $\eps_l>0$ such that $\tilde \iota(u)_\eps(K_l)\subseteq T'$
for all $\eps<\eps_l$.  Then with $\mu$ chosen with respect to this
sequence $\eps_l$ we define
$$
\iota(u)_\eps(x) := r\circ \tilde \iota(u)_{\mu(\eps,x)}(x)\,.
$$
Then clearly $\iota(u)\in \tilde \gs[X,Y]$ and $\iota$ commutes
with restrictions. Injectivity of $\iota$ follows from the fact that
$\iota(u)_\eps \to u$ uniformly on compact sets for $u\in {\mathcal
  C}(X,Y)$. Indeed, for each $l\in \N$ we may choose some $K_l'\comp
Y$ such that $\iota(u)_\eps(K_l)\cup u(K_l) \subseteq K_l'$ for $\eps$
sufficiently small. For such $\eps$ and all $x\in K_l$ we therefore
have
$$
|\iota(u)_\eps(x) - u(x)| = |r(\tilde\iota(u)_\eps(x))-r(u(x))| \le 
\|D\tilde r \|_{L^\infty(\mathrm{ch}(K_l'))} |\tilde\iota(u)_\eps(x) - u(x)|
$$
which gives the result. An analogous calculation, based on the fact
that $\tilde\iota|_{\cinfty(X,\R^s)}$ $=$ $\sigma$ establishes that
$\iota|_{\cinfty(X,Y)} = \sigma$. \ep

\brem 
At first sight it would seem that an alternative, more direct
proof of Th.\ \ref{embth1} could be carried out using Th.\ 
\ref{sheafth}: for open subsets $U$, $V$ of $\R^m$ resp.\ $\R^n$, an
embedding of $\mathcal{C}(U,V)$ into $\gs[U,V]$ can be achieved by
modifying the embedding $\iota^n$ of $\D'(U)^n$ into $\gs(U)^n$ (cf.\ 
e.g., \cite{book}, Th.\ 1.2.20) using a compact exhaustion of $U$ and
a gluing procedure as in Rem.\ \ref{mgrem}, together with the fact
that $\iota^n(u)_\eps$ converges to $u$ uniformly on compact sets. In
this way, for each $u\in \mathcal{C}(X,Y)$ one can construct a family
of embeddings of $u|_{U_\alpha}$ into $\gs[U_\alpha,Y]$ for a covering
$\{U_\alpha\mid \alpha\in A\}$ of $X$ by chart domains. Note, however,
that the resulting family in general is {\em not} coherent: in fact,
this would require the standard embedding $\iota: \mathcal{C}\to \gs$
to satisfy $f\circ \iota(u) \circ g = \iota(f\circ u\circ g)$ for
diffeomorphisms $f,\, g$ and $u$ continuous, which is manifestly wrong
in general (even for $g=\mathrm{id}$ it only holds on the level of
association, cf.\ \cite{book}, Prop.\ 1.2.70 (iv) and Sec.\ 3.2.2).
\et

As there is no notion of manifold-valued distributions (see, however,
section \ref{dprimexysec}), and in the absence of additional structure
no growth conditions can be imposed on mappings between $X$ and $Y$,
the maximal extension of $\mathcal{C}(X,Y)$ relevant to our present
considerations is the space of locally bounded measurable mappings
from $X$ to $Y$.  In order to analyze it, we shortly recall some
notions from measure theory on smooth manifolds (cf.\ 
\cite{dieudonne3}). Since zero sets are the same for any Lebesgue
measure on $X$ there is a well-defined notion of Lebesgue-measurable
subset of $X$: $A\subseteq X$ is Lebesgue-measurable iff it can be
written in the form $A=\bigcup_{j\in \N} K_j \cup N$ with $K_j\comp X$
for all $j$ and $N$ a zero set. A mapping $u:X\to Y$ is called
measurable if inverse images of Borel-measurable subsets of $Y$ under
$u$ are Lebesgue-measurable in $X$. For $X=\R^m$, $Y=\R^n$ this
precisely reproduces the usual notion of Lebesgue-measurability.
Moreover, $u$ is Lebesgue-measurable if and only if $\psi\circ u \circ
\vphi^{-1}$ is Lebesgue-measurable for any charts $\psi$ of $Y$ and
$\vphi$ of $X$.

We define ${\mathcal L}^\infty_{\mathrm{loc}}(X,Y)$ to be the set of
all Lebesgue-measurable mappings $u:X\to Y$ which are locally bounded
in the following sense: For each $K\comp X$ there exist $K'\comp Y$
and a zero set $N\subseteq X$ such that $u(K\setminus N) \subseteq
K'$. Factoring this space by the equivalence relation of coinciding
Lebesgue-almost everywhere we obtain the space
$L^\infty_{\mathrm{loc}}(X,Y)$.  In the particular case of $Y$ being a
submanifold of some $\R^s$, $L^\infty_{\mathrm{loc}}(X,Y)$ can be
identified with the subset $\tilde L_{\mathrm{loc}}^\infty(X,Y)$ of
$L^\infty_{\mathrm{loc}}(X,\R^s)$ (in the usual sense) whose elements
possess a representative mapping $X$ into $Y$. Although in the
$\R^n$-setting, $L^\infty_{\mathrm{loc}}$ can be embedded (as a
subspace of $\D'$) into $\gs$, the following example demonstrates that
the construction given in Th.\ \ref{embth1} does not carry over to
this setting in general: 

\bex\label{embex} Let $X=\R$, $Y=S^1\subseteq \R^2\cong
\C$ and let $u\in L^\infty_{\mathrm{loc}}(X,Y)$ be given by $u(x)=
(0,-1)$ for $x<0$ and $u(x)=(0,1)$ for $x>0$, i.e., $u(x)=
i\mathrm{sgn}(x)$. Then with $\rho$ a standard mollifier as above,
$u*\rho_\eps(x)$ for $x\in \R$ covers the entire line connecting $-i$
and $i$.  Therefore, it can never be contained in any tubular
neighborhood of $S^1$ and the construction breaks down.  \et

Although the above example shows that one cannot embed $u$
using the construction involving tubular neighborhoods of $S^1$ this
does not mean that one cannot construct an embedding of $L^\infty_{\mathrm{loc}}(X,S^1)$ 
into $\gs[X,S^1]$.
In \ref{embex} it is fairly clear that one can obtain such an
embedding by thinking of $S^1$ as a manifold with covering space $\R$.
We will think of points in $S^1$ as lying in the interval
$J=[0,2\pi)$. Thus an element $u \in L^\infty_{\mathrm{loc}}(X,S^1)$
defines a function $\hat u:X \to J \subset \R$ and hence an element of
$L^\infty_{\mathrm{loc}}(X,\R)$. This may be smoothed using a standard
mollifier to give a family of smooth functions $\bar u_\eps=\hat u
*\rho$. Projecting this back to $S^1$ by defining $\tilde
u_\eps=\exp(i\bar u_\eps)$ we define a family of smooth functions
whose equivalence class defines an element $\tilde u \in \gs[X,
S^1]$. Furthermore it is clear that when applied to a
smooth function this embedding gives the same result as $\sigma$. This
construction may be generalized in an obvious way to give an embedding
of $L^\infty_{\mathrm{loc}}(X,T^n)$ into $\gs[X,T^n]$.

A similar strategy may be applied to embed a function $ u \in
L^\infty_{\mathrm{loc}}(X,Y)$ where $Y$ is some compact Riemann
surface apart from $S^2$. We have excluded the 2-sphere and already
dealt with the case of the torus, so the remaining Riemann surfaces
have the structure $\Hbb^2/\Gamma$, where $\Hbb^2$ is the upper half
plane with the hyperbolic metric and $\Gamma$ is a properly
discontinuous subgroup of $PSL(2,\R)$ that acts freely on $\Hbb^2$. In
this situation there exists some fundamental polygon $F$ with finitely
many sides and for each side there is precisely one other side
obtained by the action of some element $g \in \Gamma$ with different
pairs of sides carried to each other by different elements of
$\Gamma$ (see Theorem 2.4.1 of \cite{Jost} for details). 
By adding on precisely one side from each pair to the
interior of $F$ we may obtain a region $D$ such that no two distinct
elements of $D$ are related by the action of any element $g \in
\Gamma$ and the sets $gD$ also cover $\Hbb^2$. Given $u \in 
L^\infty_{\mathrm{loc}}(X,Y)$ we may now define a corresponding function 
$\hat u: X \to D \subset \R^2$ by defining $\hat u(x)$ to be the 
unique point in $D$ such that $u(x)=g\hat u(x)$. The function $\hat u$ is an element of 
$L^\infty_{\mathrm{loc}}(X,\R^2)$ and hence may be smoothed by 
convolution with a standard mollifier to give a family of smooth
functions $\bar u_\eps=\hat u*\rho_\eps$. Projecting back down to
$Y$ gives a family of smooth functions $\tilde u_\eps$ whose
equivalence class defines an element of $\gs[X,Y]$. Again for the case
of a smooth function the result of this embedding is the same as
applying $\sigma$.

Turning now to the vector bundle setting we note that also in this situation we have
a canonical embedding $\hat\sigma: \mathrm{Hom}(E,F) \to \mathrm{Hom}_\gs(E,F)$,
$\hat\sigma: u \mapsto [(u)_\eps]$. By $\mathrm{Hom_{\mathcal C}}(E,F)$ we denote
the space of continuous vector bundle homomorphisms from $E$ to $F$. Using this notation we have:

\bt \label{vbemb} There exists an embedding $\hat\iota: \mathrm{Hom}_{\mathcal{C}}(E,F) \to \mathrm{Hom}_\gs(E,F)$ 
with the following properties:
\begin{itemize}
\item[(i)] $\hat\iota$ is a sheaf morphism.
\item[(ii)] $\hat\iota|_{\mathrm{Hom}(E,F)} = \hat\sigma$.
\item[(iii)] Any representative of $\hat\iota(u)$ converges to $u$ uniformly on compact sets. 
\end{itemize}
\et
\pr The idea is to first embed the base component by Th.\ \ref{embth1} and then employ a partition
of unity argument in the fiber components adapting the construction of an embedding of $\D'(X)$
into $\gs(X)$ given in \cite{ndg}, Th.\ 2. To this end we shall make use of the technical 
apparatus developed in Th.\ \ref{homsheaf}. Let $u\in  \mathrm{Hom}_{\mathcal{C}}(E,F)$.
With $\iota$ as in Th.\ \ref{embth1} we set $w\equiv [(w_\eps)_\eps] := \iota(\underline{u})$.
Choose now countable vector bundle atlases $\{(\Phi_i,U_i)\mid i \in \N\}$ of $E$,
$\{(\Psi_j,V_j)\mid j \in \N\}$ of $Y$ such that for each $i$ there exists
some $j$ with $\underline{u}(U_i) \subseteq V_{i_j}$. We may further suppose that each
$U_i$ is relatively compact and choose a partition of unity $\{\chi_i\mid i\in \N\}$
on $X$ with $\supp \chi_i \subset U_i$ for each $i$. Let $\zeta_i \in \D(U_i)$
such that $\zeta_i \equiv 1$ on $\supp \chi_i$ and choose some mollifier $\rho \in  
{\mathcal S}(\R^{m'})$ with unit integral and all higher moments vanishing.
Since $w_\eps$ converges to $\underline{u}$ uniformly on compact sets by Th.\ \ref{embth1},
for each $i$ there exists some $\eps_i$ such that 
\begin{equation}\label{wepssupp}
w_\eps(\supp \zeta_i)\subseteq V_{j_i}
\end{equation}
for $\eps<\eps_i$. Employing a gluing function as in the proof of Th.\ \ref{embth1} if necessary,
we may suppose without loss of generality that (\ref{wepssupp}) in fact holds for all $\eps\in (0,1]$
and all $i\in \N$.
Then define
$v_\eps^{i} \in \mathrm{Hom}_{\mathcal C}(\vphi_i((\supp\zeta_i)^\circ)\times\R^{m'},\psi_{j_i}(V_{j_i})\times \R^{n'})$ by 
$$
(x,\xi) \mapsto \left(\psi_{j_i}(w_\eps(\vphi_i^{-1}(x))), \left((\chi_i\circ\vphi_i^{-1})\cdot
(\Psi_{j_i}\circ u \circ \Phi_i^{-1})^{(2)}\right)*\rho_\eps(x)\cdot \xi\right)
$$
where convolution with $\rho_\eps$ is to be read component-wise.
Finally, we construct $v_\eps\in \mathrm{Hom}(E,F)$ by
$$
v_\eps:= \sum_{i\in \N} \zeta_i\bullet_{w_\eps} \Psi_{j_i}^{-1}\circ v_\eps^{i}\circ \Phi_i
$$
and set $\hat\iota(u):=[(v_\eps)_\eps]$. Properties (i)--(iii) then follow from Th.\ \ref{embth1} 
and the proof of \cite{ndg}, Th.\ 2. \ep

\brem It is in fact possible to embed  vector bundle homomorphisms into $\mathrm{Hom}_\gs(E,F)$ 
which are more singular (in the fiber component) than those considered in Th.\ \ref{vbemb}. To introduce such mappings
we first recall an alternative description of smooth vector bundle homomorphisms $f: E\to F$.
For vector bundle charts $\Phi_{\alpha}, \Phi_{\alpha'}$ of $E$ and $\Psi_\beta, \Psi_{\beta'}$ 
of $F$ we set $\Phi_{\alpha\alpha'}:=\Phi_\alpha\circ\Phi_{\alpha'}$ and analogously for 
$\Psi_{\beta\beta'}$. It then follows from (\ref{vbhomloc}) that 
\begin{eqnarray}
\!f_{\Psi_\beta\Phi_\alpha}^{(1)}(x)\! &=& \!\Psi_{\beta\beta'}^{(1)}(f_{\Psi_{\beta'}\Phi_{\alpha'}}^{(1)}(\Phi_{\alpha'\alpha}^{(1)}(x))) 
\label{vbtrans1}\\
\!f_{\Psi_\beta\Phi_\alpha}^{(2)}(x)\! &=& \!
\Psi_{\beta\beta'}^{(2)}(f_{\Psi_{\beta'}\Phi_{\alpha'}}^{(1)}(\Phi_{\alpha'\alpha}^{(1)}(x)))\cdot
f_{\Psi_{\beta'}\Phi_{\alpha'}}^{(2)}(\Phi_{\alpha'\alpha}^{(1)}(x))
\cdot \Phi_{\alpha'\alpha}^{(2)}(x) \label{vbtrans2}
\end{eqnarray}
so that we may identify smooth vector bundle homomorphisms $f:E\to F$ with families of
smooth local vector bundle homomorphisms $f_{\Psi_\beta\Phi_\alpha}$ which transform
according to (\ref{vbtrans1}), (\ref{vbtrans2}). As in the case of distributions on a
manifold (cf.\ \cite{hoer1}, Sec.\ 6.3) we may directly generalize the transformation
behavior (\ref{vbtrans1}), (\ref{vbtrans2}) by allowing $f_{\Psi_\beta\Phi_\alpha}^{(2)}$
to be a matrix with distributional entries and by replacing compositions with
distributional pullbacks. One restriction, however, immediately becomes apparent:
$f_{\Psi_{\beta'}\Phi_{\alpha'}}^{(1)}$ has to be supposed smooth in order for 
the right hand side of (\ref{vbtrans2}) to be well defined (i.e., to avoid ill-defined products). 
This maximal class 
of distributional vector bundle homomorphisms (with smooth base component)
can be embedded into $\mathrm{Hom}_\gs(E,F)$ by a direct adaptation of the
construction given in the proof of Th.\ \ref{vbemb}. 
\et

\section{Manifold-valued distributions} \label{dprimexysec}\label{D'(X,Y)}
Since distributions on manifolds are defined as continuous linear functionals
(on the space of compactly supported densities) there is a priori no concept
of distributions taking values in a differentiable manifold. In this section
we propose the construction of a space $\D'[X,Y]$ of distributions defined on
$X$ and taking values in the manifold $Y$. The strategy is to define 
$\D'[X,Y]$ as a quotient of a suitable subspace of $\gs[X,Y]$.

We begin by analyzing the local situation. Sequential approaches to the
theory of distributions in fact have a long history, dating back already to
\cite{Mik}, see also \cite{AMS} and \cite{Temp}. The starting point for
such considerations is the simple observation that the space of distributions
is isomorphic to the quotient of the space of nets $(u_\eps)_\eps$ of smooth
functions which converge in $\D'$ modulo the space of nets with $u_\eps \to 0$
distributionally. The Colombeau approach is of course sequential in nature and
in fact it was noted already in \cite{c2} that $\D'$ can equivalently defined as 
a certain subspace of of the Colombeau algebra $\gs$. The most direct way of
realizing $\D'$ within the Colombeau framework is the following:
\blem \label{dprimechar} Let $\Omega\subseteq\R^n$ open and set
$$
\mathcal{A}(\Omega):= \{u=[(u_\eps)_\eps]\in \gs(\Omega)\mid u_\eps \mbox{ converges in } \D'(\Omega)\}\,.
$$
Call two elements $u$, $v$ equivalent, $u\equiv v$, if $u_\eps -v_\eps\to 0$ distributionally.
Then $\D'(\Omega)$ is linearly isomorphic to $\mathcal{A}(\Omega)/\!\equiv$.
\et
\pr With $\iota: \D'(\Omega) \to \gs(\Omega)$ the standard embedding, the map
$$
\begin{array}{rcl}
\D'(\Omega) &\hookrightarrow & \mathcal{A}(\Omega)/\!\equiv \\
w & \to & [\iota(w)]_{\equiv}
\end{array}
$$
is a linear isomorphism. In fact, $\iota(w)_\eps \to w$ distributionally as $\eps\to 0$. Linearity 
and injectivity are clear. \ep

This result of course immediately generalizes to the case where the domain is a differentiable 
manifold (by employing the embedding provided by \cite{ndg}, Th. 1). However, when 
generalizing the target space to a smooth manifold $Y$, additional aspects have to be taken
into account: most importantly, diffeomorphism invariance has to be implemented. Moreover,
in the absence of additional structure it is not to be expected that unbounded distributions
can be modelled (e.g., regularizations of the delta distribution have to display
growth properties which can only be realized in the presence of scales).
As similar obstacles have already been overcome in the construction of $\gs[X,Y]$, 
in view of Lemma \ref{dprimechar} the following definition provides a natural 
generalization: 

\bd Let 
$$
\mathcal{A}[X,Y] := \{u=[(u_\eps)_\eps]\in \gs[X,Y] \mid \forall f\in \cinfty(Y),\ f\circ u_\eps \mbox{ converges  in }\ \D'(X)\}
$$
and call $u$, $v \in \mathcal{A}[X,Y]$ model-associated, $u\model v$, if $f\circ u_\eps - f\circ v_\eps \to 0$
in $\D'(X)$
for all $f\in \cinfty(Y)$. The quotient space $\D'[X,Y] := \mathcal{A}[X,Y]/\!\model$ is called the space of distributions
on $X$ taking values in $Y$.  
\et
The concept of model-association has been introduced in \cite{flows}, Sec. 5 and compared
with various other concepts of association. In particular, it was shown that $\model$ is
strictly stronger than the usual concept of association in case $Y$ is a Euclidean space.
Diffeomorphism invariance is implemented in the above definition through composition  with smooth functions $f$:
any such $f$ can be viewed as an extension of a component of a chart, hence each $f\circ u_\eps$
represents a different local picture of $u_\eps$. 

Our first observation concerning the above definition is that $u=[(u_\eps)_\eps]\in 
\mathcal{A}[X,Y]$ implies that for each $f\in \cinfty(Y)$ there exists some $u_f\in \D'(X)$
such that $f\circ u_\eps \to u_f$ in $\D'(X)$. In fact, it turns out that $u_f$ has to
be an element of $L^\infty_{\mathrm{loc}}(X)$: let $K$ be a compact subset of $X$. Then by
the c-boundedness of $(u_\eps)_\eps$ it follows that $(f\circ u_\eps)_\eps$ is
uniformly bounded on $K$, hence possesses a weak-$*$ convergent subsequence 
$(f\circ u_{\eps_k}|_K)_k$
with limit $v_f \in L^\infty(K)$. Therefore $u_f\!\!\mid_{\D(K)}=v_f$ and the claim
follows by covering $X$ with relatively compact open sets. 
Moreover, from Th.\ \ref{embth1} (iii) it follows that $\mathcal{C}(X,Y)$ can be embedded as a subspace
into $\D'[X,Y]$ and that for $u\in \mathcal{C}(X,Y)$, each $u_f$ is a continuous function.

For the special case $Y=\R$ and $f=\mathrm{id}_{\R}$ it follows that if a sequence $(u_\eps)_\eps\in\mathcal{A}[X,\R]$
then $(u_\eps)_\eps$ converges to an element $u\in L^\infty_{\mathrm{loc}}(X)$ distributionally. Moreover, if 
we denote by $\iota_S$ a sheaf embedding of $\D'(X)$ into $\G(X)$ as in \cite{book}, Thm.\ 3.2.10 we see that if
$w\in\D'(X)$, $\iota_S(w)=[(w_\eps)_\eps]$ with $(w_\eps)_\eps\in\mathcal{A}[X,\R]$ then
$w\in L^\infty_{\mathrm{loc}}(X)$. Hence representatives of elements of $\D'[X,\R]$ which come from
classical distributions actually are $L^\infty_{\mathrm{loc}}(X)$-functions.

From these considerations one might be led
to believe that $\D'[X,Y]$ singles out a certain subspace of $L^\infty_{\mathrm{loc}}(X,Y)$
in the sense that for each $u\in \mathcal{A}[X,Y]$ there should exist some underlying $v\in 
L^\infty_{\mathrm{loc}}(X,Y)$ with $u_f = f\circ v$ for all $f\in \cinfty(Y)$. However, as
the following example demonstrates, the situation is more involved:

\bex\label{4.3} 
Let $X=Y=\R$ and consider the net $u_\eps(x) = \sin(x/\eps)$. As is easily verified, 
each $u_\eps^j$ ($j\in \N_0$) converges in $\D'(\R)$. Hence every polynomial in $u_\eps$ converges
in $\D'$ and by the Weierstrass approximation theorem it follows that $(u_\eps)_\eps
\in \mathcal{A}[X,Y]$ (in fact, $f\circ u_\eps$ converges in $\D'$ for all continuous $f$). Nevertheless,
there does not exist any $v\in L^\infty_{\mathrm{loc}}(\R)$ such that $f\circ u_\eps \to f\circ v$ for
all $f\in \cinfty(\R)$. Otherwise, the choice of $f=\mathrm{id}$ would entail $v=0$,
whereas for $f(x) = x^2$, $f\circ u_\eps \to 1/2$, a contradiction.
\et 

In case $Y=\R^n$, a local description of elements of $\D'[X,Y]$ can be given in terms of 
Young measures (cf.\ e.g., \cite{Evansweak}, p.16, Th.\ 11). In fact, in this situation
we may w.l.o.g.\ (using charts) suppose that $X=U$ is a bounded open subset of $\R^m$ and that for 
$[(u_\eps)_\eps]\in \D'[U,\R^n]$, $\{u_\eps\mid \eps\in I\}$ is uniformly 
bounded (by choosing $f=\mathrm{pr}_j$). Hence there exists a subsequence $(u_{\eps_k})_k$ and a 
Young measure $(\nu_x)_{x\in U}$ such that for each $f\in \mathcal{C}(Y)$
$$
f\circ u_{\eps_k}(x) \to \int_{\R^m} f(y) \, d\nu_x(y)
$$
weak-$*$ in $L^\infty(U)$. 

Thus elements of $\D'[X,\R^n]$ do not possess an underlying description by an $L^\infty_{\mathrm{loc}}$-function
(as shown by Ex.\ \ref{4.3}) but they do possess a description in terms of Young measures in the above sense.

To conclude this section, we turn to the question of stability of $\D'[X,Y]$
under differentiation. We first note that differentiating an element of $\G[X,Y]$ gives a 
generalized vector bundle homomorphism of the respective tangent bundles rather than simply an
element of $\G[X,Y]$ (see \cite{gfvm}). Hence when considering stability under differentiation we need to take into
account this change of category. 

As a preparatory result, we provide an alternative
description of vb-equivalence for generalized vector bundle homomorphisms of tangent bundles:
\blem Let $E=TX$, $F=TY$ and $u$, $v \in {\es}_{M}^{\mathrm{VB}}[TX,TY]$. Then the following
statements are equivalent:
\begin{itemize}
\item[(i)] $u  \sim_{vb} v$.
\item[(ii)] $Tf\circ u \sim_{vb} Tf \circ v$ in 
${\es}_{M}^{\mathrm{VB}}[TX,\R\times \R^n]$ for all $f\in \cinfty(Y)$.
\end{itemize}
\et
\pr We first note that by \cite{gfvm2}, Th.\ 4.2 it suffices to prove the equivalence
for $\sim_{vb0}$ instead of $\sim_{vb}$ (i.e., with $k=0$ in Def.\ \ref{longdef} (b) (ii)).
Moreover, the same result establishes (i)$\Rightarrow$(ii) (even for general vector
bundle homomorphisms instead of tangent maps). Conversely, property (i) of 
Def.\ \ref{longdef} (b) follows from Th.\ \ref{mainchar} (ii) since for $\vphi$ a chart of
$X$ (hence $T\vphi$ a vector bundle chart of $TX$) we have
$$
(Tf\circ u_\eps - Tf\circ v_\eps)^{(1)}_{\mathrm{id}T\vphi}(\vphi(p))
= (f\circ \underline{u_\eps} - f\circ \underline{v_\eps})\circ\vphi^{-1}(p)\,.
$$
For establishing property (ii) it suffices to note that any vector bundle chart
$T\psi$ of $TY$ is the restriction of some $Tf$ for a suitable extension of
$\psi$ (cf.\ the proof of Prop.\ 4.1 in \cite{gfvm2}).\ep

The above result suggests that to test equivalence of tangent maps of elements of
$\D'[X,Y]$ one should compose with tangent maps of elements of $\cinfty(Y)$.
This leads to the following
\bd Two elements $u$, $v$ of $\mathrm{Hom}_\gs(TX,TY)$ are called model-vb-equivalent, 
if $Tf\circ u_\eps - Tf\circ v_\eps$ converges to $0$ distributionally for all 
$f\in \cinfty(Y)$.
\et
More precisely, this means that for any chart $(\vphi,V)$ of $X$ and any $f\in \cinfty(Y)$
we have 
\beas
&&f\circ u_\eps\circ \vphi^{-1} - f\circ v_\eps\circ \vphi^{-1} \to 0 \\
&& D(f\circ u_\eps\circ \vphi^{-1} - f\circ v_\eps\circ \vphi^{-1}) \to 0
\eeas
in $\D'(\vphi(V))$. Since distributional convergence is stable under derivatives
it therefore follows that for $u$, $v \in \mathcal{A}[X,Y]$, model equivalence
entails vb-model equivalence of the respective tangent maps. As in the local
distributional setting, stability under differentiation can therefore be achieved
also for $\D'[X,Y]$ if iterated tangent maps of elements of $\cinfty(Y)$ are
used for testing equivalence of derivatives. 

{\em Acknowledgments:} We would like to thank James D.E. Grant for several helpful discussions.


\begin{thebibliography}{10}

\bibitem{AMS}
{Antosik, P., Mikusinski, J., Sikorski, R.}
\newblock {\em Theory of Distributions -- The Sequential Approach}.
\newblock Elsevier, 1973.

\bibitem{AB}
{Aragona, J., Biagioni, H.~A.}
\newblock Intrinsic definition of the {C}olombeau algebra of generalized
  functions.
\newblock {\em Analysis Mathematica}, {\bf 17}:75--132, 1991.

\bibitem{Biag}
{Biagioni, H.~A.}
\newblock {\em A Nonlinear Theory of Generalized Functions}, volume~{\bf 1421}
  of {\em Lecture Notes in Mathematics}.
\newblock Springer, Berlin, 1990.

\bibitem{c1}
{Colombeau, J.~F.}
\newblock {\em New Generalized Functions and Multiplication of Distributions}.
\newblock North Holland, Amsterdam, 1984.

\bibitem{c2}
{Colombeau, J.~F.}
\newblock {\em Elementary Introduction to New Generalized Functions}.
\newblock North Holland, Amsterdam, 1985.

\bibitem{Cbull}
{Colombeau, J.~F.}
Multiplication of distributions.
{\em Bull. Amer. Math. Soc. (N.S.)}, {\bf 23}:251--268, 1990.

\bibitem{CM}
{Colombeau, J.~F., Meril, A.}
\newblock Generalized functions and multiplication of distributions on
  {${\mathcal C}^\infty$} manifolds.
\newblock {\em J.~Math.~Anal.~Appl.}, {\bf 186}:357--364, 1994.

\bibitem{deR}
{De Rham, G.}
\newblock {\em Differentiable Manifolds}, volume~{\bf 266} of {\em Grundlehren
  der mathe\-mati\-schen Wissenschaften}.
\newblock Springer, Berlin, 1984.

\bibitem{RD}
{De Roever, J.~W., Damsma, M.}
\newblock Colombeau algebras on a {${\cal C}^\infty$}-manifold.
\newblock {\em {Indag.~Mathem., N.S.}}, {\bf 2}(3), 1991.

\bibitem{dieudonne3}
{Dieudonn\'e, J.}
\newblock {\em Treatise on Analysis}, volume~{\bf 3}.
\newblock Academic Press, New York, 1972.

\bibitem{DP} Dapi\'c, N.; Pilipovi\'c, S. Microlocal analysis of Colombeau's 
generalized functions on a manifold.  {\em Indag. Math. (N.S.)}  7  no. 3, 293--309, 1996.
 
\bibitem{DP2} Dapi\'c, N.; Pilipovi\'c, S.; Scarpalezos, D. Microlocal analysis of 
Colombeau's generalized functions: propagation of singularities.  {\em J. Anal. Math.}
75  51--66, 1998.

\bibitem{DHPV}
Delcroix, A., Hasler, M., Pilipovi\'c, S., Valmorin, V. 
Generalized function algebras as sequence space algebras.  
{\em Proc. Amer. Math. Soc.}  132  (2004),  no. 7, 2031--2038

\bibitem{Evansweak}
{Evans, L. C.}
\newblock {\em Weak convergence methods for nonlinear partial differential
  equations}, volume~74 of {\em CBMS Regional Conference Series in
  Mathematics}.
\newblock 1990.

\bibitem{found}
{Grosser, M., Farkas, E., Kunzinger, M., Steinbauer, R.}
\newblock On the foundations of nonlinear generalized functions {I}, {II}.
\newblock {\em Mem. Amer. Math. Soc.}, {\bf 153}(729), 2001.

\bibitem{book}
{Grosser, M., Kunzinger, M., Oberguggenberger, M., Steinbauer, R.}
\newblock {\em Geometric Theory of Generalized Functions}, volume 537 of {\em
  Mathematics and its Applications 537}.
\newblock Kluwer Academic Publishers, Dordrecht, 2001.

\bibitem{vim}
{Grosser, M., Kunzinger, M., Steinbauer, R., Vickers, J.}
\newblock A global theory of algebras of generalized functions.
\newblock {\em Adv. Math.}, 166:179--206, 2002.

\bibitem{hirsch}
{Hirsch, M. W.}
\newblock {\em Differential topology}, volume~33 of {\em Graduate Texts in
  Mathematics}.
\newblock Springer-Verlag, New York, 1994.

\bibitem{hoer1}
{H\"ormander, L.}
\newblock {\em The Analysis of Linear Partial Differential Operators {I}},
  volume~{\bf 256} of {\em Grundlehren der mathematischen Wissenschaften}.
\newblock Springer, Berlin, 1990.

\bibitem{Jel}
{Jel\'\i nek, J.}
\newblock An intrinsic definition of the {C}olombeau generalized functions.
\newblock {\em Comment. Math. Univ. Carolinae}, {\bf 40}:71--95, 1999.

\bibitem{Jel2} Jel\'\i nek, J. On introduction of two diffeomorphism invariant Colombeau algebras.  
{\em Comment. Math. Univ. Carolin.}  45  (2004),  no. 4, 615--632.

\bibitem{Jel3} Jel\'\i nek, J. Equality of two diffeomorphism invariant Colombeau algebras.  
{\em Comment. Math. Univ. Carolin.}  45  (2004),  no. 4, 633--662.

\bibitem{Jost} Jost, J.
 Compact Riemann surfaces. Springer, Berlin, 1997.

\bibitem{gfvm}
{Kunzinger, M.}
\newblock Generalized functions valued in a smooth manifold.
\newblock {\em Monatsh.\ Math.}, {\bf 137}:31--49, 2002.

\bibitem{symm}
{Kunzinger, M., Oberguggenberger, M.}
\newblock Group analysis of differential equations and generalized functions.
\newblock {\em SIAM J. Math. Anal.}, {\bf 31}(6):1192--1213, 2000.

\bibitem{flows}
{Kunzinger, M., Oberguggenberger, M., Steinbauer, R., Vickers, J.}
\newblock Generalized flows and singular {ODE}s on differentiable manifolds.
\newblock {\em Acta Appl. Math.}, {\bf 80}:221--241, 2004.

\bibitem{ndg}
{Kunzinger, M., Steinbauer, R.}
\newblock Foundations of a nonlinear distributional geometry.
\newblock {\em Acta Appl. Math.}, 71:179--206, 2002.

\bibitem{gprg}
{Kunzinger, M., Steinbauer, R.}
\newblock Generalized {pseudo-}{R}iemannian geometry.
\newblock {\em Trans. Amer. Math. Soc.}, 354(10):4179--4199, 2002.

\bibitem{gfvm2}
{Kunzinger, M., Steinbauer, R., Vickers, J.}
\newblock Intrinsic characterization of manifold-valued generalized functions.
\newblock {\em Proc.\ London Math.\ Soc.}, {\bf 87}(2):451--470, 2003.

\bibitem{lee}
{Lee, J. M.}
\newblock {\em Introduction to smooth manifolds}, volume 218 of {\em Graduate
  Texts in Mathematics}.
\newblock Springer-Verlag, New York, 2003.

\bibitem{marsden}
{Marsden, J.~E.}
\newblock Generalized {H}amiltonian mechanics.
\newblock {\em Arch.\ Rat.\ Mech.\ Anal.}, {\bf 28}(4):323--361, 1968.

\bibitem{marti} Marti, Jean-Andre $({\mathcal C},{\mathcal E},{\mathcal P})$-sheaf 
structures and applications.  
Nonlinear theory of generalized functions (Vienna, 1997),  175--186, Chapman \& Hall, Boca Raton, FL, 1999.

\bibitem{Mik}
{Mikusi\'nski, J.}
\newblock Sur la m\'ethode de g\'en\'eralisation de {L}aurent {S}chwartz sur la
  convergence faible.
\newblock {\em Fund.~Math.}, {\bf 35}:235--239, 1948.

\bibitem{N}
{Nash, J.}  The imbedding problem for Riemannian manifolds.
 Ann. of Math. (2)  {\bf 63}  (1956), 20--63.

\bibitem{MOBook}
{Oberguggenberger, M.}
\newblock {\em Multiplication of Distributions and Applications to Partial
  Differential Equations}, volume~{\bf 259} of {\em Pitman Research Notes in
  Mathematics}.
\newblock Longman, Harlow, U.K., 1992.

\bibitem{Schw}
{Schwartz, L.}
\newblock Sur l'impossibilit\'e de la multiplication des distributions.
\newblock {\em C.~R.~Acad.~Sci.~Paris}, {\bf 239}:847--848, 1954.

\bibitem{SV}
{Steinbauer, R., Vickers, J.}
\newblock The use of generalized functions and distributions in general
  relativity. {\em Classical Quantum Gravity}  23  (2006),  no. 10, R91--R114.

\bibitem{Temp}
{Temple, G.}
\newblock Theories and applications of generalized functions.
\newblock {\em J.~London Math.~Soc.}, {\bf 28}:134--148, 1953.

\end{thebibliography}
\end{document}